\documentclass[12pt,reqno]{amsart}
\usepackage[cp1251]{inputenc}
\usepackage[english]{babel}
\usepackage{amsfonts}
\usepackage{amssymb,amsmath,amsthm,eufrak,euscript}
\usepackage{graphicx}
\usepackage[usenames]{color}
\usepackage{colortbl}
\newtheorem{lemma}{Lemma}
\newtheorem{theorem}{Theorem}

\newtheorem{remark}{Remark}

\numberwithin{equation}{section}

\begin{document}

\sloppy

\begin{center}
{\bf \Large Finite-gap potentials and integrable geodesic equations on a 2-surface\footnote{Both authors are supported by the grant of the Russian Science Foundation No. 24-11-00281, https://rscf.ru/project/24-11-00281/}}
\end{center}

\medskip

\begin{center}
{\bf S.V. Agapov, A.E. Mironov}
\end{center}

\medskip

\begin{quote}
\noindent{\bf Abstract. }{\small In this paper we show that the one-dimensional Schr\"odinger equation can be viewed as the geodesic equation of a certain metric on a 2-surface. In case of the Schr\"odinger equation with a finite-gap potential, the metric and geodesics are explicitly found in terms of the Baker--Akhiezer function.
}

\medskip


\end{quote}

\hfill {\it To the $70$th anniversary of S.V. Bolotin and the $60$th anniversary of D.V. Treschev}

\medskip

\section{Introduction and the main result}

{\it Geodesics} of a Riemannian (a pseudo-Riemannian) metric $ds^2=g_{ij}(u)du^idu^j$ on a 2-surface $M$ in local coordinates $u^1,$ $u^2$ satisfy a system of equations of the following form:
\begin{equation}
\label{1}
\ddot{u}^k + \Gamma^k_{ij}\dot{u}^i\dot{u}^j = 0, \qquad k=1,2,
\end{equation}
where $\Gamma^k_{ij}$ are the Christoffel symbols.
The search for Riemannian metrics on 2-surfaces such that the equations~\eqref{1} can be integrated in quadratures is a classical and, generally speaking, very complicated problem. It has long been known that an integrability takes place only in exceptional cases. For instance, in the famous work~\cite{1} Jacobi integrated the geodesic equations on the three-axis ellipsoid in terms of elliptic functions. In the overwhelming majority of cases, the geodesic equations are not integrable (e.g., see~\cite{2}).

Let us rewrite~\eqref{1} in the form of a Hamiltonian system
\begin{equation}
\label{2}
\dot{u}^k=\frac{\partial H}{\partial p_k}, \qquad \dot{p}_k=-\frac{\partial H}{\partial u^k}, \qquad k=1,2
\end{equation}
with the Hamiltonian $H = g^{ij}(u)p_ip_j / 2;$ projections of solutions to the system~\eqref{2} on $M$ coincide with the geodesics, that is, they satisfy~\eqref{1}. Recall that any function $F(u, p)$ which is preserved along trajectories of the system~\eqref{2}, i.e.
$$
\frac{dF}{dt} = \sum_{i=1}^2 \left( \frac{\partial F}{\partial u^i} \frac{\partial H}{\partial p_i} - \frac{\partial F}{\partial p_i} \frac{\partial H}{\partial u^i} \right) \equiv 0
$$
is called its {\it first integral}.
Notice that the Hamiltonian $H$ is the first integral of~\eqref{2}. Therefore, due to the Liouville---Arnold theorem (\cite{3}), the equations~\eqref{2} (and, consequently,~\eqref{1}) can be integrated in quadratures if there exists an additional first integral $F(u, p)$ which is functionally independent on $H$ almost everywhere. A large number of papers is devoted to the search for such integrals, mainly polynomial in momenta (e.g., see~\cite{4}---\cite{6} and references within).

Notice that the question of integrability of other Hamiltonian systems is also relevant. Let us mention some known general results. For instance, the so-called natural systems are closely related to geodesic flows; the relation is given by the well-known Maupertuis's principle. Topological obstacles to existence of analytical first integrals of such systems on 2-surfaces were described in~\cite{Kz1} (see also~\cite{Bol1}). In~\cite{Tm1} this result was generalized to the case of multidimensional manifolds. In~\cite{KzTr1} natural systems on multidimensional tori were investigated; a criterion for existence of a complete set of polynomial integrals was found in the case of trigonometric potentials.

\begin{remark}
The Liouville---Arnold theorem indicates only the principal possibility to integrate the Hamilton's equations having a complete set of first integrals. Explicit integration of these equations is an another problem which often turns out to be not easier than a direct search for the first integrals.
\end{remark}

Denote for convenience $u^1=x,$ $u^2=y.$ We shall consider $y$ as a function of $x$ along a geodesic. Then the system~\eqref{1} implies that
\begin{equation}
\label{3}
\frac{d^2y}{dx^2} = A_3(x,y) \left( \frac{dy}{dx} \right)^3 + A_2(x,y) \left( \frac{dy}{dx} \right)^2 + A_1(x,y) \frac{dy}{dx} + A_0(x,y),
\end{equation}
where coefficients $A_j(x,y),$ $j=0, \ldots, 3$ are expressed in terms of the Christoffel symbols $\Gamma^i_{jk}$ of the initial metric in the following way:
\begin{equation}
\label{4}
A_3=\Gamma^1_{22}, \qquad A_2=2\Gamma^1_{12}-\Gamma^2_{22}, \qquad A_1=\Gamma^1_{11}-2\Gamma^2_{12}, \qquad A_0=-\Gamma^2_{11}.
\end{equation}

Now consider the equation~\eqref{3} with arbitrary coefficients $A_j(x,y),$ $j=0, \ldots, 3.$ We shall call it {\it metrisable} if the relations~\eqref{4} hold true, i.e. if there exists a metric such that its geodesics satisfy~\eqref{3}. A criterion of metrisability of equations of the form~\eqref{3} was obtained in~\cite{7}. It is the following. Consider an over-determined system of linear PDEs on the functions $\psi_j(x,y)$, $j=1, \ldots, 3:$
\begin{equation}
\label{5}
 \begin{gathered}
(\psi_{1})_x=\frac23 A_1\psi_1-2A_0\psi_2,\hfill\\
(\psi_3)_y=2 A_3\psi_2-\frac23A_2\psi_3,\hfill\\
(\psi_{1})_y+2(\psi_{2})_x=\frac43 A_2\psi_1-\frac23 A_1\psi_2-2A_0\psi_3,\hfill\\
(\psi_{3})_x+2(\psi_{2})_y=2 A_3\psi_1-\frac43 A_1\psi_3+\frac23A_2\psi_2.\hfill
\end{gathered}
\end{equation}
If this system has a solution such that $\Delta=\psi_1\psi_3-\psi_2^2\not\equiv0$, then the equation~\eqref{3} defines the geodesics of the metric
$$
 ds^2=\frac{1}{\Delta^2}(\psi_1dx^2+2\psi_2dxdy+\psi_3dy^2)
 $$
 and the relations~\eqref{4} hold true.

Many well-known equations of mathematical physics have the form ~\eqref{3} with appropriate coefficients $A_j(x,y).$ The question about metrisability of these equations is interesting but, as far as we know, has been little studied so far. Certain results on metrisability of similar equations and their first integrals can be found in~\cite{8},~\cite{9}.

In this paper we study the metrisability problem for the Schr\"odinger equation
\begin{equation}
\label{6}
 L_2y=\frac{d^2y}{dx^2} - u(x) y = z y,\quad y=y(x).
\end{equation}
We shall integrate the over-determined system of PDEs~\eqref{5} and thereby find the metric
$$
 ds^2=g_{11}(x,y)dx^2+2g_{12}(x,y)dxdy+g_{22}(x,y)dy^2
$$
on a surface in case when the geodesic equation~\eqref{3} takes the form~\eqref{6}. The following theorem holds true.
\begin{theorem}
The Schr\"odinger equation~\eqref{6} is metrisable, and coefficients of the metric $ds^2$ have the form
\begin{equation}
	\label{g1}
	g_{11}(x,y)=\lambda(x,y) ( r_0 + \frac{1}{2} y ( yl''(x)-4s'(x)-2yl(x)(u(x)+z))),
\end{equation}
\begin{equation}
	\label{g2}
	g_{12}(x,y)=\lambda(x,y) (s(x)-\frac{1}{2}yl'(x)),\quad g_{22}(x,y)= \lambda(x,y) l(x),
\end{equation}
where
$$
\lambda(x,y) = \left( \left( s(x)-yl'(x)/2 \right)^2-l(x) \left( r_0+y (yl''(x)-4s'(x)-2yl(x)(u(x)+z))/2 \right) \right)^{-2},
$$
$r_0$ is a constant. Here functions $s(x), l(x)$ satisfy the equations
\begin{equation}
	\label{8}
	\frac{d^2s(x)}{dx^2}-u(x) s(x)=zs(x),\quad
\end{equation}
\begin{equation}
	\label{9}
	\frac{d^3l(x)}{dx^3} - 4(u(x)+z)\frac{d l(x)}{dx}-2u'(x)l(x)=0.
\end{equation}
The Gaussian curvature $K$ of the metric $ds^2$ has the form
$$
4 K = 2 \left( r_0l-s^2 \right) l'' - r_0(l')^2+ 4sl's'-4l(s')^2 + 4l \left( s^2-r_0l \right)(u+z).
$$
\end{theorem}

The equation~\eqref{9} is well-known in the theory of the one-dimensional Schr\"odinger operator $L_2$. It appears in various situations related to $L_2$ (e.g., see~\cite{GD},~\cite{LP}, an equation~\eqref{00} below).
Let
$$
 s(x)=a_1s_1(x)+a_2s_2(x),\quad a_1,a_2\in\mathbb{R}
$$
be a general solution to the equation~\eqref{8}, where $s_1(x),s_2(x)$ are certain solutions. Then one can check that a general solution to the equation~\eqref{9} has the form
$$
 l(x)=b_1s_1^2(x)+b_2s_1(x)s_2(x)+b_3s_2^2(x), \quad b_j\in\mathbb{R}.
$$
This is one of the beautiful properties of the equation~\eqref{9}. Thus the metric $ds^2$ depends on six parameters $r_0,a_i,b_j$. Moreover, for any values of parameters the metrics are geodesically equivalent, that is, non-parameterized geodesics coincide for any parameters. If $r_0=a_1=a_2=0,$ then the Gaussian curvature $K=0$. Geodesically equivalent metrics were studied, for instance, in~\cite{Di},~\cite{MT}.

In case of finite-gap operators solutions to the equations~\eqref{8},~\eqref{9} as well as metric coefficients are expressed in terms of the Baker--Akhiezer function.

\begin{remark}
There exists a construction which connects geodesics on multidimensional quadrics with the one-dimensional finite-gap Schr\"odinger operator and is based on other ideas (e.g., see~\cite{Neumann}---\cite{VVV}). More precisely, it is proved in~\cite{Kn} that after an appropriate parametrization geodesics on quadrics can be transformed by the Gauss map into trajectories of the classical Neumann problem on the unit sphere with a quadratic potential; moreover, any solution to the Neumann problem can be obtained in this way by choosing an appropriate quadric. It is proved in~\cite{VVV} that any one-dimensional Schr\"odinger operator with a finite-gap potential can be constructed via certain solution to this Neumann problem.
\end{remark}

In Section 2 we recall some results about one-dimensional finite-gap Schr\"odinger operators which we need further. In particular, we recall a construction of the Baker--Akhiezer function which is needed to describe geodesics. In Section 3 we give the proof of Theorem 1 and consider the Lame operator as an example.

\section{Finite-gap potentials and the Baker--Akhiezer function}

\subsection{The Baker--Akhiezer function}

Let us recall the definition of a finite-gap operator. Firstly suppose that a potential $u(x)$ of the operator
$L_2=\frac{d^2}{dx^2} - u(x)$ is real, smooth and periodic $u(x+\tau)=u(x), \tau$ is a period.
An eigenfunction $\varphi(x)$ of the operator $L_2$ is called the Bloch one
$$
L_2\varphi=z\varphi,\quad z\in\mathbb{R}
$$
if
$$
\varphi(x+\tau)=e^{ip}\varphi(x),\quad p\in\mathbb{R}.
$$
A Floquet spectrum of the operator $L_2$ consists of eigenvalues $z$ corresponding to the Bloch functions. In general case the Floquet spectrum of the operator $L_2$ with a periodic potential $u(x)$ is a union of an infinite number of segments. The potential $u(x)$ (the operator $L_2$) is called {\it finite-gap} if its Floquet spectrum consists of a finite number of segments, one of which is semi-infinite. The most important results in the theory of one-dimensional finite-gap Schr\"odinger operators were obtained in~\cite{IM}--\cite{D}, in particular, in~\cite{IM} the potential $u(x)$ was found in terms of the theta function of the Jacobian variety of the hyperelliptic spectral curve (see below). Novikov \cite{N} proved that if a one-dimensional smooth periodic Schr\"odinger operator $L_2$ commutes with a certain differential operator of order $2g+1$ of the form
$$
L_{2g+1}=\partial_x^{2g+1}+u_{2g-1}(x)\partial_x^{2g-1}+\dots+u_0(x),
$$
then the operator $L_2$ is the finite-gap one. The innverse statement was proved by Dubrovin \cite{D}, namely, if $L_2$ is finite-gap, then $L_2$ commutes with a certain operator $L_{2g+1}$. In general case, if $L_2$ commutes with $L_{2g+1}$, then the potential $u(x)$ is quasi-periodic. In what follows, by the finite-gap operator $L_2$ (correspondingly the potential $u(x)$) we will mean a real operator that commutes with a certain operator of an odd order.

If the operator $L_2$ commutes with $L_{2g+1}$, then due to Burchnal---Chaundy lemma~\cite{BC} operators $L_2, L_{2g+1}$ satisfy the equation $F(L_2,L_{2g+1})=0$, where
$$
 F(z,w)=w^2-(z^{2g+1}+c_{2g}z^{2g}+\dots+c_0),
$$
$c_j$ are some constants. A spectral curve $\Gamma$ of a pair of operators $L_2, L_{2g+1}$ is defined as a smooth compactification with a point at infinity $q=\infty$ of an affine curve given in $\mathbb{C}^2$ by the equation
\begin{equation}
\label{sc}
 w^2=z^{2g+1}+c_{2g}z^{2g}+\dots+c_0.
\end{equation}
The spectral curve parameterizes common eigenfunctions of the operators, namely, if
$$
 L_2\psi=z\psi,\quad L_{2g+1}\psi=w\psi,
$$
then the point $P=(z,w)$ lies on the spectral curve, i.e.
$$
 \psi=\psi(x,z,w)=\psi(x,P),\quad P\in\Gamma.
$$
Common eigenfunctions $\psi(x,P)$ are {\it the Baker--Akhiezer functions} which were found in~\cite{K} in general case for ordinary commuting differential operators of rank one.
Such functions are constructed via the spectral data
$$
\{ \Gamma, q, k^{-1}, \gamma_1, \ldots, \gamma_g \},
$$
where $\Gamma$ is a Riemannian surface of genus $g,$ $q \in \Gamma$ is a marked point, $k^{-1}$ is a local parameter in a neighborhood of the point $q,$ $\gamma_1+\ldots+\gamma_g$ is a non-special divisor. There exists the only one function $\psi(x,P),$ $P\in \Gamma,$ which satisfies the following conditions:

1) in the neighborhood of $q$ the function $\psi$ has the form
$$
\psi = e^{xk} \left( 1+\frac{\xi_1(x)}{k}+\frac{\xi_2(x)}{k^2}+\ldots \right),
$$

2) the function $\psi$ is meromorphic on $\Gamma \setminus \{q\}$ with simple poles in $\gamma_1, \ldots, \gamma_g.$

${\ }$

\noindent {\bf Example 1.}
Let $\Gamma$ be an elliptic curve
$$
\Gamma = \mathbb{C}/\Lambda,\quad \Lambda=\{2mw_1+2nw_2, \ n, m \in \mathbb{Z}\}.
$$
the Baker--Akhiezer function related to the spectral data $\{ \Gamma, 0, -z, -\gamma\}$ has the form
\begin{equation}\label{BA}
\psi(x,z) = e^{-x\zeta(z)} \frac{\sigma(z+x+\gamma)\sigma(\gamma)}{\sigma(z+\gamma)\sigma(x+\gamma)},
\end{equation}
and
\begin{equation}\label{Lame}
\left( \partial^2_x-2\wp(x+\gamma) \right) \psi(x,z) = \wp(z) \psi(x,z).
\end{equation}
Here $\zeta,$ $\sigma,$ $\wp$ are the Weierstrass elliptic functions. The operator $\partial^2_x-2\wp(x+\gamma)$ commutes with the operator $$L_3=\partial^3_x-3\wp(x+\gamma)\partial_x-\frac{3}{2}\wp'(x+\gamma).$$

In the general case the Baker--Akhiezer function is expressed in terms of the theta function of the Jacobian variety of the spectral curve $\Gamma$ \cite{K}. Let us choose a basis of cycles $a_1,\dots,a_g,b_1,\dots,b_g$ on $\Gamma$ with the following intersection numbers:
$$
 a_i\circ a_j=b_i\circ b_j=0,\quad a_i\circ b_j=\delta_{ij}.
$$
Let $\omega_1,\dots,\omega_g$ be the basis of Abelian differentials on $\Gamma$ normalized by the condition
$$
 \int_{a_i}\omega_j=\delta_{ij}.
$$
We denote by $\Omega$ the matrix of $b$-periods with components $\Omega_{ij}=\int_{b_i}\omega_j$. The matrix
$\Omega$ is symmetric $\Omega_{ij}=\Omega_{ji}$ and has a positive-definite imaginary part
${\rm Im} \Omega>0$. The Jacobian variety of the Riemann surface $\Gamma$ has the form $J(\Gamma)={\mathbb C}^g/\{\mathbb{Z}^g+\Omega\mathbb{Z}^g\}$. The theta function of the Jacobian variety is given by the series
$$
 \theta(z)=\sum_{n\in\mathbb{Z}^g}\exp(\pi i n\Omega n^t+2\pi i n z^t),\quad  z=(z_1,\dots,z_g)\in\mathbb{C}^g.
$$
The theta function has the following property
$$
 \theta(z+m\Omega +n)=\exp(-\pi i m\Omega m^t-2\pi i m z^t)\theta(z),\quad n,m\in{\mathbb Z}^g.
$$
The Abel map $A:\Gamma\rightarrow J(\Gamma)$ is given by the formula
$$
 A(P)=\left(\int_{P_0}^P\omega_1,\dots,\int_{P_0}^P\omega_g\right),
$$
where $P_0\in\Gamma$ is a certain fixed point. We denote by $\Omega_1$ a meromorphic 1-form on $\Gamma$ with the only pole of order two at $q$ and zero $a$-periods $\int_{a_i}\Omega_1=0$, denote by $U$ the vector of its $b$-periods
$$
 U=\left(\int_{b_1}\Omega_1,\dots,\int_{b_g}\Omega_1\right).
$$
The Baker--Akhiezer function has the form
$\psi=\frac{\tilde{\psi}}{\xi_0(x)}$, where
$$
 \tilde{\psi}=\exp(2\pi ix\int_{P_0}^P\Omega_1)\frac{\theta(A(P)-A(\gamma)-K_{\Gamma}+xU)}{\theta(A(P)-A(\gamma)-K_{\Gamma})},
$$
here $A(\gamma)=A(\gamma_1)+\dots+A(\gamma_g)$, $K_{\Gamma}$ is the vector of Riemann constants, $\xi_0(x)$ is the expansion coefficient of the function $\tilde{\psi}(x)$ in a neighborhood of $q$
$$
 \tilde{\psi} = e^{xk} \left( \xi_0(x)+\frac{\xi_1(x)}{k}+\frac{\xi_2(x)}{k^2}+\ldots \right).
$$
The finite-gap potential has the form
$$
 u(x) = 2\partial^2_x \log \theta (xU+V) + c, \qquad V=-A(\gamma)+K_{\Gamma},
$$
$c$ is a constant.

Smooth real potentials $u(x)$ are distinguished as follows. We assume that all branch points of the spectral curve $\Gamma$ are real, i.e. the equation (\ref{sc}) takes the form
$$
 w^2=(z-z_1)\dots(z-z_{2g+1}),\quad  z_j\in\mathbb{R},
$$
then $c_j\in\mathbb{R}$ and $\Gamma$ admits an anti-holomorphic involution
$$
 \tau: \Gamma\rightarrow\Gamma,\quad \tau(z,w)=(\bar{z},-\bar{w}).
$$
In this case, $\Gamma$ is an $M$-curve, i.e. it has $g$ maximal possible cycles that are stationary with respect to $\tau$. Assume that each of the points $\gamma_1,\dots,\gamma_g$ lies on its own stationary cycle. Then the potential $u(x)$ is smooth and real. Conditions on the spectral data for periodic operators are found in \cite{MO} (see also \cite{GS}).
Note that in case of a periodic potential, periodic and antiperiodic solutions to the equation (\ref{6}) correspond to the branch points $R_j=(z_j,0)$ of the spectral curve. Also note that if $z\ne z_j$, then the basis in the space of solutions to the equation $L_2y=zy$ has the form $\psi(x,z,w), \psi(x,z,-w)$. In case of a periodic potential and periodic solutions $L_2y=z_jy$ the basis has the form
\begin{equation}
	\label{P}
 \psi(x,z_j,0),\quad \partial_{k_j}\psi(x,z_j,0),
\end{equation}
where $k_j$ is a local parameter in a neighborhood of the branch point $R_j$, $k_j=\sqrt{z-z_j}$.

\subsection{Finite-gap Schr\"odinger operators}
Let $L_2$ be a finite-gap operator. Then the operator $L_2-z$ is factorized in the following way (\cite{15}):
$$
L_2-z = -\left( \partial_x+\chi_0(x,P) \right) \left( \partial_x-\chi_0(x,P) \right), \quad P\in\Gamma,
$$
where $\chi_0(x, P) = \frac{\psi_x(x,P)}{\psi(x,P)}$ is a meromorphic function on $\Gamma$ with simple poles at some points
$$
 P_1(x)=(z_1(x),w_1(x)), \ldots, P_g(x)=(z_g(x),w_g(x))
$$
and $q$. Points $P_j(x)$ are zeros of the Baker--Akhiezer function. The function $\chi_0$ has the form
$$
 \chi_0 = \frac{Q_x}{2Q}+\frac{w}{Q}, \quad Q(x,z)=z^g+\alpha_{g-1}(x)z^{g-1}+\ldots+\alpha_0(x),
$$
where $\alpha_{j}(x)$ are some functions.
The polynomial $Q$ in $z$ satisfies the equation
$$
4w^2=4(z+u(x))Q^2+Q_x^2-2QQ_{xx}.
$$
After differentiating this equality with respect to $x$ and dividing by $Q$ we obtain
\begin{equation}
\label{00}
Q_{xxx}-4(z+u(x))Q_x-2u_x(x)Q=0.
\end{equation}
We note that $z_j(x)$ are zeros of $Q$, i.e.
\begin{equation}
	\label{Q}
 Q=(z-z_1(x))\dots(z-z_g(x)).
\end{equation}
In case when $P$ coincides with the branch point $R_j$ we have
$$
 \chi_0(x,P_j)=\frac{Q_x(x,z_j)}{2Q(x,z_j)}=\frac{\psi_x(x,P_j)}{\psi(x,P_j)}.
$$
This implies that
$$
 Q(x,z_j)=\delta\psi^2(x,R_j),
$$
where $\delta$ is a constant.

${\ }$

\noindent {\bf Example 2.}
One example of the finite-gap operator is given by the Lame operator
$$
 \partial_x^2-g(g+1)\wp(x+\omega_2), \qquad g \in \mathbb{N},
$$
where $\wp(x)$ is the Weierstrass elliptic function with the lattice of periods $\Lambda=\{2mw_1+2nw_2, \ n, m \in \mathbb{Z}\}.$
If $\omega_1\in\mathbb{R}$ and $\bar{\Lambda}=\Lambda$, then the Lame operator is real, smooth and periodic.
The function $\wp(x)$ satisfies the equation of the form
\begin{equation}\label{wp}
 \left( \wp'(x) \right)^2 = 4 \wp^3 - g_2 \wp(x) -g_3.
\end{equation}
The function $Q$ for the Lame operator is a polynomial of degree $g$ in $\wp(x+\omega_2)$ as well as in $z$ (\cite{M1}, \cite{M2})
$$
Q = B_g\wp^g(x+\omega_2)+B_{g-1}(z)\wp^{g-1}(x+\omega_2)+\ldots+B_0(z),
$$
where coefficients $B_s$ can be found by the following recurrence formula
$$
B_s=\frac{(s+1) \left( 8B_{s+1}z-B_{s+2}g_2(s+2)(2s+3)-2B_{s+3}g_3(s+2)(s+3) \right)}{4(2s+1)(g^2+g-s(s+1))},
$$
here $B_g$ is a constant, $B_m=0$ for all $m>g.$
The spectral curve of the Lame operator is given by the following equation
$$
w^2=\frac{1}{4} \left( 4B_0^2z-B_0(4B_2g_3+B_1g_2)+B_1^2g_3 \right).
$$
Another way to construct the spectral curve of the Lame operator was proposed in~\cite{GV}.

\section{Proof of Theorem 1}

Let us prove Theorem 1. First of all, we prove the following lemma.

\begin{lemma}
Assume that
\begin{equation}
\label{7}
A_0(x,y) = (u(x)+z) y, \qquad A_1(x,y)=A_2(x,y)=A_3(x,y)=0.
\end{equation}
Then a solution to the system~\eqref{5} has the following form:
$$
\psi_3(x,y)=l(x), \qquad \psi_2(x,y) = s(x)-\frac{1}{2}yl'(x),
$$
$$
\psi_1(x,y) = r_0 + \frac{1}{2} y \left( yl''(x)-4s'(x)-2yl(x)(u(x)+z) \right),
$$
where $r_0$ is an arbitrary constant and functions $s(x),$ $l(x)$ satisfy the equations~\eqref{8}, \eqref{9}.
\end{lemma}

Let us prove Lemma 1. Consider the system~\eqref{5} provided that the coefficients $A_0, \ldots, A_3$ have the form~\eqref{7}. Integrating the second, the fourth and the third equations~\eqref{5} successively, we obtain
$$
\psi_3(x,y)=l(x), \qquad \psi_2(x,y) = s(x)-\frac{y}{2}l'(x),
$$
$$
\psi_1(x,y) = r(x) + \frac{1}{2} y \left( yl''(x)-4s'(x)-2yl(x)(u(x)+z) \right),
$$
where $l(x),$ $s(x),$ $r(x)$ are arbitrary functions. After that the first equation of the system~\eqref{5} takes the form
$$
 r'(x)-2y(s''(x)-(u(x)+z)s(x))+\frac{1}{2}y^2(l'''(x)-4(u(x)+z)l'(x)-2u'(x)l(x))=0.
$$
This implies that $r(x) = r_0$ is an arbitrary constant and functions $s(x),$ $l(x)$ satisfy the equations~\eqref{8},~\eqref{9}. Lemma 1 is proved.

To complete the proof of Theorem 1 it is left to notice that for the functions $\psi_1, \psi_2, \psi_3$ found in Lemma 1 the metric coefficients
$$
g_{11}(x,y) = \frac{\psi_1}{(\psi_1 \psi_3-\psi_2^2)^2}, \qquad g_{12}(x,y) = \frac{\psi_2}{(\psi_1 \psi_3-\psi_2^2)^2}, \qquad g_{22}(x,y) = \frac{\psi_3}{(\psi_1 \psi_3-\psi_2^2)^2},
$$
have the form \eqref{g1}, \eqref{g2}. Theorem 1 is proved.

In case of finite-gap operators the equation (\ref{9}) coincides with the equation (\ref{00}). As a solution to the equation (\ref{9}) one can take a linear combination of the functions $\psi^2(x,z,w)$, $\psi^2(x,z,-w)$ and $Q(x,z)$ (see \eqref{Q}).
In case of the Lame operator the function $Q(x,z)$ is defined in Example 2. As a solution to the equation (\ref{8}) we take the linear combination
$$
 s(x)=\beta_1\psi(x,z,w)+\beta_2\psi(x,z,-w),\quad \beta_1, \beta_2\in\mathbb{R};
$$
and in case of periodic solutions ($z=z_j$) we take a linear combination of the functions~\eqref{P}.
Thus, all the functions $\psi_1(x,y),$ $\psi_2(x,y),$ $\psi_3(x,y)$ are found.

To sum it up, the geodesics of the metric~\eqref{g1}, \eqref{g2} satisfy the Schr\"odinger equation~\eqref{6} which in case of a finite-gap periodic potential $u(x)$ is integrable in terms of the Baker---Akhiezer function. In addition notice that the curves $x=const$ are also geodesics since, as can be easily verified, $\Gamma^1_{22} \equiv 0.$

${\ }$

\noindent {\bf Example 3.} Let us find the metric $g_{ij}(x,y)$ in case of the Lame operator for $g=1$, $u(x)=2\wp(x+\gamma)$.
In this case one may take the function $s(x)$ in the form
$$
 s(x)= e^{-x\zeta(z)} \frac{\sigma(z+x+\gamma)\sigma(\gamma)}{\sigma(z+\gamma)\sigma(x+\gamma)}.
$$
The equation~\eqref{9} takes the form
$$
 \frac{d^3l(x)}{dx^3} - 4(2\wp(x+\gamma)+\wp(z))\frac{d l(x)}{dx}-4\wp'(x+\gamma)l(x)=0.
$$
This equation has the following solution
$$
 l(x)=\wp(z)-\wp(x+\gamma).
$$
Now the metric is given by the formulae~\eqref{g1},~\eqref{g2}.

Let us consider the degeneration of the elliptic spectral curve into a rational curve with a cuspidal special point. For that, we put $g_2=g_3=0$ for the Weierstrass function in the equation \eqref{wp}. Then the Weierstrass elliptic functions degenerate into elementary ones:
$$
 \hat{\wp}(x) = \frac{1}{x^2}, \qquad \hat{\zeta}(x) = \frac{1}{x}, \qquad \hat{\sigma}(x) = x.
$$
The Baker--Akhiezer function~\eqref{BA} transforms into the function
$$
 \hat{\psi}(x,z)=e^{-x/z} \frac{ (x+z+\gamma)\gamma}{(x+\gamma) (z+\gamma)},
$$
and the Lame operator takes the form
$$
L_2= \partial^2_x-\frac{2}{(x+\gamma)^2}.
$$
A general solution to the equation
$$
L_2 \varphi = \frac{1}{z^2} \varphi
$$
has the form
$$
\varphi = \beta_1 e^{-x/z} \frac{ (x+z+\gamma)}{(x+\gamma) (z+\gamma)} + \beta_2 e^{x/z} \frac{ (x-z+\gamma)}{(x+\gamma) (-z+\gamma)}, \quad \beta_1, \beta_2 \in \mathbb{R}.
$$
For simplicity we put $\gamma=0, z=1$, then the metric takes the form
$$
g_{11} = \frac{r_0x^4-y^2(x^4+x^2+1)}{(x^2y^2-r_0x^2+r_0)^2}, \quad
g_{12} = -\frac{xy}{(x^2y^2-r_0x^2+r_0)^2}, \quad
g_{22} = \frac{x^2(x^2-1)}{(x^2y^2-r_0x^2+r_0)^2},
$$
its Gaussian curvature is constant and equal to
$$
 K=-r_0.
$$
We have
$$
\varphi = \beta_1 e^{-x} \frac{1+x}{x} + \beta_2 e^{x} \frac{1-x}{x}.
$$
Certain examples of geodesics corresponding to various values of parameters $\beta_1,$ $\beta_2$ are given below on Fig. 1 --- 3.

\begin{figure}[h]
\begin{center}
\begin{minipage}[h]{0.31\linewidth}
\includegraphics[width=0.9\linewidth]{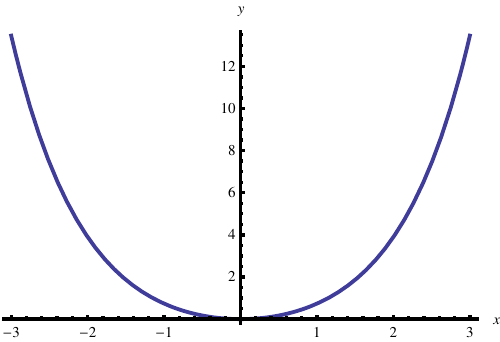}
\caption{$\beta_1=1,$ $\beta_2=-1;$}
\label{ris:experimoriginal}
\end{minipage}
\hfill
\begin{minipage}[h]{0.31\linewidth}
\includegraphics[width=0.9\linewidth]{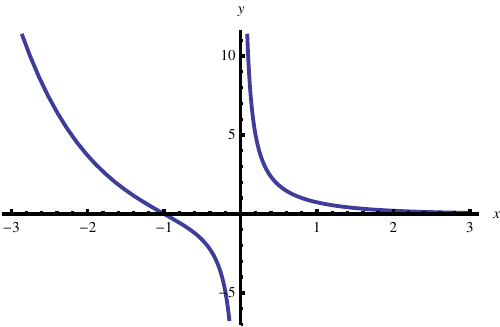}
\caption{$\beta_1=1,$ $\beta_2=0;$}
\label{ris:experimcoded}
\end{minipage}
\hfill
\begin{minipage}[h]{0.31\linewidth}
\includegraphics[width=0.9\linewidth]{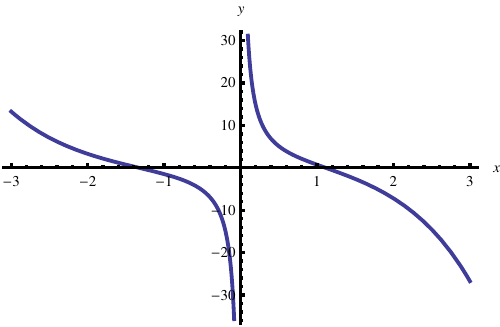}
\caption{$\beta_1=1,$ $\beta_2=2$}
\label{ris:experimcoded}
\end{minipage}
\hfill
\end{center}
\end{figure}

\vspace{0.2cm}
{\bf Acknowledgments.} The authors thank an anonymous referee for pointing out to us the connection between geodesic flows on quadrics and one-dimensional finite-gap Schr\"odinger operators described in Remark 2.

${\ }$

Sergei Agapov,

Novosibirsk State University, 630090, 1, Pirogova str., Novosibirsk;

Sobolev Institute of Mathematics SB RAS,  630090, 4 Acad. Koptyug avenue

e-mail: agapov.sergey.v@gmail.com

${\ }$

Andrey Mironov,

Novosibirsk State University, 630090, 1, Pirogova str., Novosibirsk;

Sobolev Institute of Mathematics SB RAS,  630090, 4 Acad. Koptyug avenue

e-mail: mironov@math.nsc.ru

\end{document}